\documentclass[10pt,reqno]{amsart}

\usepackage{amsmath,amssymb,amsthm}
\usepackage{graphicx}
\usepackage[all]{xypic}
\input xypic
\newtheorem{problem}{Problem}[section]
\newtheorem{example}{Example}[section]
\newtheorem{remark}{Remark}
\newtheorem{theorem}{Theorem}
\newtheorem{corollary}{Corollary}
\newtheorem{lemma}{Lemma}                                                                                                                                                                                      \newtheorem{Proposition}{Proposition}[section]
\newtheorem{definition}{Definition}
\newtheorem{proposition}{Proposition}
\numberwithin{equation}{section}
\def\nga{\widetilde}
\def\ex{\begin{example}
  }
  \def\eex{\end{example}}
\def\thr{\begin{theorem}}
\def\ethr{\end{theorem}}
\def\pro{\begin{Proposition}}
\def\epro{\end{Proposition}}
\def\coro{\begin{corollary}}
\def\ecoro{\end{corollary}}
\def\df{\begin{definition}}
\def\edf{\end{definition}}
\def\lm{\begin{lemma}}
\def\elm{\end{lemma}}
\def\pf{\begin{proof}}
\def\epf{\end{proof}}
\def\problem{\begin{problem}}
\def\eproblem{\end{problem}}

\def\it{\begin{itemize}}
\def\hit{\end{itemize}}

\newcommand{\vc}{\infty}
\newcommand{\vm}{\forall}
\newcommand{\tru}{\setminus}
\newcommand{\td}{\Longleftrightarrow}
\newcommand{\sr}{\longrightarrow}

\newcommand{\ps}{\dfrac}

\newcommand{\mtn}{\rightarrow}
\newcommand{\gd}{\underline}
\newcommand{\gt}{\overline}

\newcommand{\seq}[1]{\left<#1\right>}

\def\Def{\mathcal{D}ef}
\def\id {\mbox{id}}
\def\isodef{\underline{\mathcal{D}ef}}
\def\Der{\mbox{Der}}
\def\Hom{\mbox{Hom}}
\def\Ext{\mbox{Ext}}
\def\Coker{\mbox{Coker}}
\def\ker{\mbox{ker}}
\def\Im{\mbox{Im}}
\def\Ann{\mbox{Ann}}

\def\c{\Bbb C}

\def\k{\Bbb K}

\def\hoa{\mathcal}
\def\xhoa{\mathcal{X}}
\def\yhoa{\mathcal{Y}}
\def\ohoa{\mathcal{O}}

\begin{document}
\title[On deformations of maps and curve singularities]
{On deformations of maps and curve singularities}

\author{Gert-Martin Greuel}

\address{Fachbereich Mathematik, Universit\"{a}t Kaiserslautern, Germany}

\email{greuel@mathematik.uni-kl.de}

\author{Cong Trinh Le}

\address{Fachbereich Mathematik, Universit\"{a}t Kaiserslautern, Germany}

\email{cle@mathematik.uni-kl.de}


\begin{abstract}
We study several deformation functors associated to the normalization
of a reduced curve singularity $(X,0) \subset (\c^n,0)$. The main new
results are explicit formulas, in terms of classical invariants of
$(X,0)$, for the cotangent cohomology groups $T^i, i = 0,1,2,$ of
these functors. Thus we obtain precise statements about smoothness
and  dimension of the corresponding local moduli spaces. We apply the
results to obtain explicit formulas resp. estimates for the
$\hoa{A}_e$-codimension of a parametrized curve singularity, where
$\hoa{A}_e$ denotes the Mather-Wall group of left-right equivalence.
\end{abstract}
\maketitle

\section{Introduction}
Deformations of singularities play an important role in several
branches of mathematics, in particular in complex analysis and in
algebraic and arithmetic geometry. The main general results of
deformation theory are the existences of a semi-universal deformation
for isolated singularities and  an infinitesimal theory of
deformations and obstructions. Beside this, not much more is known
except we have more information, as it is, for example, the case for
complete intersections or for Cohen-Macaulay singularities in
codimension two.

The present paper treats deformations of reduced curve singularities
where we have, as additional information, the  parametrization. That
is, if $X \subset (\c^n,0)=:Y$ is the germ of a reduced curve
singularity, we have a commutative diagram
 $$\xymatrix@C=12pt@R=8pt@M=4pt{
\gt{X}\ar[d]_{n}  \ar[dr]^\varphi &  \\
X~ \ar @{^{(}->}[r]_>>>>{j}~ & ~ Y=(\c^n,0).} $$ Here $n$ is the
normalization map, $\gt{X}:=(\gt{X},\gt{0}),$ with
$\gt{0}=n^{-1}(0),$ is a multigerm which can be identified with a
multigerm $(\c,S)$ where $S \subset \c$ is finite set of cardinality
equal to the number of branches of $X$, $j$ is the inclusion map and
the composition $\varphi = j\circ n$ is a (primitive) parametrization
of $X$.

Now we cannot only deform $X$ but also the (multi-)germs $\gt{X}$ and
$Y$ and the morphisms $n, \varphi$ and $j$. The relation between
these deformations has been  elaborated by Buchweitz in his thesis
(cf. \cite{Bu}, not published) and is partially reproduced in the
textbook \cite{GLS} to which we refer as the main reference. Since
$\gt{X}$ and $Y$ are smooth, the deformation theory of $\varphi$ is
rather simple and we can use the relation between the different
deformation functors to obtain useful information about deformations
of $X$. This was done in \cite[Chapter II.2]{GLS} for plane curve
singularities where the infinitesimal deformation spaces $T^1$ and
the obstruction spaces $T^2$ for the various deformations have been
computed and where the results were applied to equinormalizable and
equisingular deformations.

Plane curve singularities are special in the sense that every
deformation of $\varphi$ induces a deformation of $X$ (cf.
\cite[Prop. 2.23]{GLS}) which is not true if $X$ is not plane, the
case which we treat in the present paper. Therefore we cannot expect
the same intimate relation between deformations of the
parametrization $\varphi$ and the deformations of $X$. However, we
show that it is nevertheless possible to compute, resp. estimate, the
dimensions of the various infinitesimal deformation and obstruction
spaces in terms of more familiar invariants of the curve singularity.
The general procedure is very similar to the plane curve case, but
the proofs are  technically more involved. Although our main results
in Theorem \ref{md1} concern the vector spaces $T^1$ and $T^2$, we
can easily deduce results about the dimension of the base spaces of
the corresponding semi-universal deformations (cf. Remark \ref{rm3.2}
at the end of Section 4).

The results of Theorem \ref{md1} are applied to obtain exact formulas
resp. estimates for the $\hoa{A}_e$-codimension of an analytic
morphism $\varphi: (\c,S) \mtn (\c^n,0)$, where $\hoa{A}_e$ denotes
the Mather-Wall group of (extended) left-right equivalence. This has
been studied before by several authors in varying generality,
including the recent paper \cite{HRR}, which has been one of the
motivations to publish our more general results. We show first that
the $\hoa{A}_e$-codimension of $\varphi$ is finite if and only if
$\varphi$ is a primitive parametrization of a reduced curve
singularity and, in this case, we can compute resp. estimate the
$\hoa{A}_e$-codimension by showing that it is equal to $\dim_\c
T^1_{(\c,S)\mtn (\c^n,0)}$ and then apply Theorem \ref{md1}. By
specifying the general result to plane curve resp. to Gorenstein
resp. to quasihomogeneous singularities etc., we recover and
generalize all (as far as we are aware of) previously known results
about $\hoa{A}_e$-codimensions.

Although we prefer a geometric point of view and work over $\c$, all
results hold (with the same  proofs) for algebroid curve
singularities over an algebraically closed field $\k$  of
characteristic $0$. In this case the deformations in Section 2 have
to be defined dually for diagrams of (complete) analytic local rings.

In positive characteristic the results do not hold as stated but have
to be modified (see \cite{CGL} for the plane curve case).

\section{Deformation functors of maps}

We state some basic facts about deformations of maps and fix
notations.  We refer to the text book \cite[Chapter II, Appendix
C]{GLS} as the main reference where further references can be found.

We consider the category of (pointed) complex spaces or complex
(multi) germs, the distinguished point is usually denoted by $0$.
$T_\epsilon$ denotes  the non-reduced point defined by
$\c[\epsilon]/(\epsilon^2)$.

Let $f: X \sr Y$ be a morphism of spaces and  $S = (S,0)$  a pointed
space or a germ. A deformation $F: \mathcal{X} \sr \mathcal{Y} $ of
$f : X \sr Y$ over $S$ is a Cartesian diagram

 \hspace{4cm}$\xymatrix@C=8pt@R=7.5pt@M=4.5pt{
\ar @{} [dr] |{\Box}
{X~~} \ar[d]_f \ar @{^{(}->}[r]& ~{\mathcal{X}}\ar[d]^F\\
\ar @{} [dr] |{\Box}
{Y~~} \ar[d] \ar @{^{(}->}[r]& ~\mathcal{Y}\ar[d]\\
\{0\} ~~ \ar @{^{(}->}[r] & ~S} $\\
such that $\hoa{X}$ and $\hoa{Y}$ are flat over $S$.  If $F': \xhoa'
\sr \yhoa' $ is another deformation over $S$, a morphism from $F$ to
$F'$ is given by morphisms $\Phi: \hoa{X} \sr \hoa{X}'$ and $\Psi:
\hoa{Y} \sr \hoa{Y}'$ over $S$ such that the obvious diagram
commutes. If $\Phi$ and $\Psi$ are isomorphisms, we have an
isomorphism of deformations over $S$. Moreover, for a morphism
$\varphi: (T,0) \sr (S,0)$, the pullback of $F$ defines a deformation
$\varphi^*F: \varphi^*\xhoa \sr \varphi^*\yhoa \sr T$ of $f$ over
$(T,0)$  which is called \emph{the induced deformation}.

We consider the following deformation categories (all deformations
are deformations over $S$):\\
\begin{tabular}{ll}
 $\Def_{X \mtn Y}(S):$ & deformations of $f : X \sr
Y;$ as defined above\\
& (X,Y  and f are deformed);\\
$\Def_{X/Y}(S): $ & deformations of $X$ over $Y$; consisting of
diagrams\\
& with $\yhoa = Y \times S$   and morphisms with
 $\Psi = \id_{Y\times S} $\\
 & ($X$ and $f$, but not $Y$, are deformed);\\
 $\Def_{ X }(S)$: & deformations of $X$; same as $\Def_{X / Y}$ but
 with  $Y$\\
 & $=\{0\}$(forgetting the middle part of the diagram);\\
 $\Def_{X \tru Y}(S)$: & deformations of $Y$ under $X$, consisting of
 diagrams \\
&  with $\xhoa = X \times S $  and morphisms with $\Phi = \id_{X\times S}$ \\
 & ($Y$ and $f$, but not $X$, are deformed);\\
$\Def_{Y}(S)$: & deformations of $Y$, same as $\Def_{X}$\\
& (forgetting the upper part of the diagram);\\
$\Def_{X\tru X \mtn Y/Y}(S)$: & deformations of $f$, also called
unfoldings,  with $\xhoa$\\
&$= X\times S,$  $ \yhoa = Y\times S$  and morphisms with $\Phi=$\\
& $\id_{X\times S},$ $ \Psi=\id_{Y\times S}$ \\
& ($f$, but neither $X$ nor $Y$, are deformed).\\
\end{tabular}

The set of isomorphism classes of $\Def_{X \mtn Y}(S)$ is denoted by
$\underline{\Def}_{X \mtn Y} (S)$ and $S \mapsto \gd{\Def}_{X\mtn
Y}(S) $ defines a functor from the category of pointed complex spaces
or complex space germs to the category of sets. It is called the
\emph{deformation functor of $X \sr Y.$ } Similar notations are used
for the other functors.

Instead of one morphism, we may deform a finite sequence $X_0 \mtn
X_1 \mtn \cdots \mtn X_k$ of morphisms or diagrams of morphisms. Such
deformations over $S$ are given by a commutative diagram with
Cartesian squares for each morphism in the diagram as above and
spaces  flat over $S$. Morphisms of such deformations are defined by
an obvious generalization.

In particular, a deformation of $\{0\} \sr X \sr Y$ over $S$ is the
same as a deformation $\xhoa \sr \yhoa$ of $ X \sr Y$ over $S$,
together with compatible sections $S \sr \yhoa$ and $ S \sr \xhoa$.
They are, therefore, also called \emph{deformations of $X \sr Y$ with
sections.} For details see \cite[Section  II. 2.3]{GLS}.

The deformation functors defined above are related in the following
diagram:
$$\xymatrix@C=14pt@R=8pt@M=3pt {
& & \isodef_{X \tru X \mtn Y/Y} \ar[dl] \ar[dr] & & \\
\isodef_Y & \isodef_{X\tru Y} \ar[l] \ar[dr] & & \isodef_{X/Y}
\ar[dl] \ar[r] & \isodef_{X}.\\
& & \isodef_{X \mtn Y} & & }$$

The natural transformations between these functors are obviously
defined by forgetting, or imposing, conditions.

\section{Cotangent braids for maps and  curve singularities}
For each of the above functors $\isodef$ there exists a cotangent
complex $L^{\bullet}$, whose $i-$th cohomology group is denoted by
$T^i$, where
$T^0, T^1, T^2$ have special meanings (cf. \cite[Appendix C]{GLS}):\\
$T^0:\quad $ vector space of compatible derivations;\\
$T^1: \quad$ $\isodef(T_\epsilon)$, vector space of infinitesimal
deformations;\\
$T^2: \quad $ vector space which contains the obstructions against
lifting to higher order. \\
 Note that isomorphic deformation functors
have isomorphic $T^0$ and $T^1$, but not necessarily isomorphic
$T^2$.

Actually, $T^i$ is a functor which has as argument a module. For
example:\\
 if $X$ is reduced we have $$T^0_X(M)\cong
 \Hom_{\ohoa_X}(\Omega_X^1, M) \mbox{ and } T^1_X(M) \cong
 \Ext_{\ohoa_X}^1(\Omega_X^1,M),$$
if $X$ is normal we have $T^2_X(M) \cong
\Ext_{\ohoa_X}^2(\Omega_X^1,M)$ and,  if $X$ is non-singular,
$T^i_X(M) = 0$ for $i\geq 1$, for any $\ohoa_X$-module $M$
($\Omega_X^1$ denotes the module of K\"{a}hler differentials). The
argument $M$ is omitted if $M$ is  the structure sheaf $\ohoa_X$.

Buchweitz \cite{Bu} observed that $T^i$ of the six functors defined
above are related nicely by a braid of four long exact sequences
$\xymatrix{ {} \ar@{=>}[r] & {}}$,$ \xymatrix{ {} \ar@{:>}[r] & {}}$,
$\xymatrix{ {} \ar@{-->}[r] & {}}$, $\xymatrix{ {} \ar@{->}[r] & {}}
$ (cf. \cite[p.446]{GLS}). Moreover,
\begin{equation}\label{pt1}
 T^i_{X\tru X \mtn Y/Y}
\cong T^{i-1}_{Y}(f_*\ohoa_X), \end{equation} for $i\geq 0$ ( $= 0$
for $i = 0$).  For proofs  see \cite[2.4]{Bu}.

We have another isomorphisms
\begin{equation}\label{pt2}
 T^i_{{X}\tru Y} \cong T^{i-1}_Y\big(f_*\ohoa_X/\ohoa_Y\big),
\end{equation}
for $i\geq 0$,    under the assumption that the structure map
$f^\sharp: \ohoa_Y \sr f_*\ohoa_X$ is injective. To prove
(\ref{pt2}), we refer \cite[Lemma 2.29, p.312]{GLS} the proof for
reduced plane curve singularities, noting that the argument is same
for the general situation. In particular, if $n: \gt{X} \sr X $ is
the normalization map of a reduced space $X$, $n^\sharp: \ohoa_X \sr
\ohoa_{\gt{X}}$ is injective, hence we have
\begin{equation}\label{pt2'}
T^i_{\gt{X}\tru X} \cong T^{i-1}_X\big(\ohoa_{\gt{X}}/\ohoa_X\big).
\end{equation}
For the deformation theoretic meaning of $T^1 $ and $T^2$ see Remark
\ref{rm3.2}.

Consider now the diagram associated to a reduced curve singularity
 $$\xymatrix@C=15pt@R=11pt@M=6pt{
(\gt{X},\gt{0})=:\gt{X}\ar[d]_{n}  \ar[dr]^\varphi &  \\
(X,0)=:X~ \ar @{^{(}->}[r]_>>>>{j}~ & ~ Y=(\c^n,0).} $$

We are going to use the following notations:\\
\begin{tabular}{ll}
$X$: & the germ of a reduced curve singularity with analytic
branches\\
& $X_1,\cdots, X_r$;\\
$\gt{X}:$ & the multigerm  $(\gt{X},\gt{0})$  of the normalization of
$X$, the
disjoint union \\
& of $\gt{X_1}, \cdots, \gt{X_r}$, the normalizations of $X_1,
\cdots,
X_r$, respectively;\\
$n$: & the normalization map of $X$ with $\gt{0} = n^{-1}(0), r = \sharp \gt{0}$;\\
$j$: & the embedding of $X$ in $Y = (\c^n,0)$;\\
$\varphi$: & the parametrization of $X$;\\
$\ohoa_Y$: & $\cong \c\{ x_1,\cdots,x_n\} =: \c\{\underline{x}\}, $
the local ring of $Y$;\\
$\ohoa_X$:& $\cong \c\{\underline{x}\}/I$, the local ring of $X$
(with $I$ the ideal defining $X$);\\
 $\ohoa_{\gt{X}}$: & $ \cong
\bigoplus\limits_{i=1}^r \c\{t_i\},$ the semi-local ring of
$\gt{X}$;\\
$\varphi^{\sharp}:$ & $\ohoa_Y \sr \ohoa_{\gt{X}}$, the
$\c-$algebra map of $\varphi$;\\
$n^{\sharp}:$ & $\ohoa_X \sr \ohoa_{\gt{X}}$, the
$\c-$algebra map of $n$.\\
\end{tabular}

We identify $\ohoa_Y$ with $\c\{\gd{x}\}$ and $\ohoa_{\gt{X}}$ with
$\bigoplus\limits_{i=1}^r \c\{t_i\}$ and  write the structure
morphism $\ohoa_Y \sr \ohoa_{\gt{X}} $ explicitly as
$$\varphi^\sharp = (\varphi^{\sharp}_i)_{i=1}^r : \c\{\gd{x}\} \sr
\bigoplus\limits_{i=1}^r \c\{t_i\}, \quad \varphi^\sharp(x_j) =:
{\varphi^{(j)}} =: (\varphi^{(j)}_i)_{i=1}^r, $$
$$\varphi^\sharp_i(x_j) = \varphi^{(j)}_i(t_i), \quad i = 1 ,\cdots, r, \quad  j =
1,\cdots,n.$$ Note that $\ker(\varphi^\sharp) = I, $ and that
$n^\sharp: \c\{\gd{x}\}/I \sr \bigoplus\limits_{i=1}^r \c\{t_i\}$ is
injective.

Let us  further denote by \\
\begin{tabular}{ll}
$\delta = \dim_\c \ohoa_{\gt{X}}/\ohoa_X$: & the $\delta$-invariant
of $X$;\\
$\tau =\dim_\c T^1_X$:& the Tjurina number of $X$;\\
$\mu = 2\delta - r + 1$: & the Milnor number of $X$;\\
$mt$:& the multiplicity of $X$;\\
$mt_i$: & the multiplicity of $X_i$;\\
 $e = 3 \delta - m_1$: & the Deligne number of $X$, where $m_1 = \dim_\c T^0_{\gt{X}}/T^0_X$.\\
\end{tabular}\\
Recall that $e$ is the dimension of a smoothing component of the base
space of the semi-universal deformation of $X$ if $X$ is smoothable
(cf. \cite{De}, \cite{Gr}).

In order to give an explicit description of the cotangent cohomology
for the parametrization $\varphi: \gt{X}\sr Y $ and the normalization
$n: \gt{X} \sr X$, and compute $T^1$, $T^2$, we also need $T^0$,
which we describe first:

\lm \label{bd1}  With the above notations, we have

\begin{enumerate}
\item[1.] $T^0_{\gt{X}\mtn Y} = \{(\xi,\eta)\in
\Der_\c(\ohoa_{\gt{X}},\ohoa_{\gt{X}}) \times
\Der_\c(\ohoa_{Y},\ohoa_{Y})| \xi\circ \varphi^{\sharp}=
\varphi^{\sharp}\circ \eta\}$,\\
$T^0_{\gt{X}\mtn X} \cong T^0_X$.
\item[2.] $T^0_{\gt{X}/X}\cong T^0_{\gt{X}/Y} = 0.$
\item[3.] $T^0_{\gt{X}\tru \gt{X} \mtn X/X} \cong T^0_{\gt{X}\tru \gt{X} \mtn
Y/Y} = 0.$
\item[4.] $T^0_{\gt{X}\tru Y}=\{\eta\in \Der_\c(\ohoa_Y,
\ohoa_Y)|\varphi^\sharp\circ \eta = 0\} \cong
\bigoplus\limits^n_{j=1}
\ohoa_Y I\tfrac{\partial }{\partial x_j},$\\
$T^0_{\gt{X}\tru X} = 0.$
\item[5.]  $T^0_X \cong \Der_\c(\ohoa_X,\ohoa_X) \cong
\Hom_\c(\Omega_X^1,\ohoa_X).$
\item[6.] For each $\ohoa_{\gt{X}}$-module  $N$, respectively each
$\ohoa_Y$-module $M$, we have
 $$T^0_{\gt{X}}(N) \cong
\bigoplus\limits_{i=1}^r N\tfrac{\partial}{\partial t_i}, \quad
T^0_Y(M) \cong \bigoplus\limits_{j=1}^n M\tfrac{\partial}{\partial
x_j}.$$ Moreover, $T^0_{\gt{X}} \cong T^0_{\gt{X}}(\ohoa_{\gt{X}}),
\quad T^0_Y \cong  T^0_Y(\ohoa_Y).$
\end{enumerate}
\elm

\pf 1. is proved in the same way for plane curve singularities as
\cite[Lemma 2.28 (1), p.309]{GLS}. Note that for $T^0_{\gt{X}\mtn X}
\cong T^0_X$
we use that $\mbox{char}(\c) = 0$.  \\
2. Since $n$ and $\varphi$ are finite maps, unramified outside $\{0\}
\in Y$, we have that $\Omega_{\gt{X}/X}^1$ and $\Omega_{\gt{X}/Y}^1$
are torsion modules and hence
$$T^0_{\gt{X}/X} = \Hom_{\ohoa_{\gt{X}} }(\Omega_{\gt{X}/X}^1,
\ohoa_{\gt{X}}) = 0,$$
$$
T^0_{\gt{X}/Y} =  \Hom_{\ohoa_{\gt{X}} }(\Omega_{\gt{X}/Y}^1,
\ohoa_{\gt{X}}) = 0.$$
 3. follows just from the definition, or  from the isomorphisms (\ref{pt1}).\\
 4.  By definition, we have $T^0_{\gt{X}\tru Y}=\{\eta\in \Der_\c(\ohoa_Y,
\ohoa_Y)|\varphi^\sharp\circ \eta = 0\}$. Since $\ker(\varphi^\sharp)
= \ohoa_Y I$, we obtain the result for $T^0_{\gt{X}\tru Y}$.
$T^0_{\gt{X}\tru X} = 0$ follows in the same way since $n^\sharp$ is
injective.\\
5. is just from the definition while 6. follows since
$\Omega^1_{\gt{X}} $ and $\Omega^1_Y$ are free due to the smoothness
of $\gt{X}$ and $Y$. \epf

 Now we consider the braids for $n : \gt{X}
\sr X$ and for $\varphi: \gt{X} \sr Y$. Since $\gt{X}$ and $Y$ are
nonsingular, we have $T^i_{\gt{X}}(-) = T^i_{Y}(-) = 0, $ for $i \geq
1$. By using  Lemma \ref{bd1} and this fact,  we get the cotangent
braids for the parametrization and the normalization as shown in
Figure \ref{dig1} and Figure \ref{dig2}.

\begin{figure} \hspace{1cm}$\xymatrix@C=1pt@R=2pt@M=6pt{
& 0 \ar@{=>}@/_4pc/@<-2ex>[dd] & & & & & & &\\
& & & &  0 \ar@{-->}[dlll]\ar@{->}[drrr]& & & & & \\
& T^0_{\gt{X}\tru Y} \ar@{-->}@/_4pc/@<-2ex>[dd]\ar @{=>}[drrr] & & &
& & & 0 \quad  \ar @{:>}[dlll]\ar@{->}@/^4pc/@<-2ex>[dd]& \\
& & & &  T^0_{\gt{X} \mtn Y} \ar@{:>}[dlll]\ar@{=>}[drrr]& & & & & \\
& T^0_{ Y} \ar@{:>}@/_4pc/@<-2ex>[dd]\ar
@{-->}[drrr]^{\varphi^{\star}} & & &
& & & T^0_{\gt{X}} \ar @{->}[dlll]^{\dot{\varphi}}\quad \ar@{=>}@/^4pc/@<-2ex>[dd]& \\
& & & &  T^1_{\gt{X}\tru \gt{X} \mtn Y/Y} \ar@{->}[dlll]\ar@{-->}[drrr]& & & & & \\
& T^1_{ \gt{X}/Y} \ar@{->}@/_4pc/@<-2ex>[dd]\ar @{:>}[drrr] & & &
& & & T^1_{\gt{X}\tru Y}\quad  \ar @{=>}[dlll]\ar@{-->}@/^4pc/@<-2ex>[dd]& \\
& & & &  T^1_{\gt{X}\mtn Y} \ar@{=>}[dlll]\ar@{:>}[drrr]& & & & & \\
& 0 & & & & & & \quad 0 \quad &
 }$
 \caption{The cotangent braid for the parametrization $\varphi:\gt{X} \mtn Y$
 }\label{dig1}

\hspace{1cm}$\xymatrix@C=1pt@R=2pt@M=6pt{ &&&&&&&&\\
&&&&&&&&\\
 & 0
\ar@{-->}@/_4pc/@<-2ex>[dd]\ar @{=>}[drrr] & & &
& & & 0 \quad \ar @{:>}[dlll]\ar@{->}@/^4pc/@<-2ex>[dd]& \\
& & & &  T^0_{\gt{X}\mtn X} \ar@{:>}[dlll]_{\cong}\ar@{=>}[drrr]& & & & & \\
& T^0_{X} \ar@{:>}@/_4pc/@<-2ex>[dd]\ar @{-->}[drrr]^{n^{\star}} & &
&& & & T^0_{\gt{X}}\quad  \ar @{->}[dlll]_{\dot{n}}\ar@{=>}@/^4pc/@<-2ex>[dd]& \\
& & & &  T^1_{\gt{X}\tru\gt{X} \mtn X/X} \ar@{->}[dlll]\ar@{-->}[drrr]& & & & & \\
&T^1_{\gt{X}/X}\quad \ar@{->}@/_4pc/@<-2ex>[dd]\ar@{:>}[drrr]
& & && & & T^1_{\gt{X}\tru X} \quad  \ar @{=>}[dlll]\ar@{-->}@/^4pc/@<-2ex>[dd]& \\
& & & &  T^1_{\gt{X}\mtn X} \ar@{=>}[dlll]\ar@{:>}[drrr]& & & & & \\
& 0 \ar@{=>}@/_4pc/@<-2ex>[dd]\ar @{->}[drrr] & & &
& & & T^1_{X}\quad  \ar @{-->}[dlll]\ar@{:>}@/^4pc/@<-2ex>[dd]& \\
& & & &  T^2_{\gt{X}\tru\gt{X} \mtn X/X} \ar@{-->}[dlll]\ar@{->}[drrr]& & & & & \\
& T^2_{\gt{X}\tru X}\ar@{-->}@/_4pc/@<-2ex>[dd]\ar @{=>}[drrr] & & &
& & & T^2_{\gt{X}/X} \quad \ar @{:>}[dlll]\ar@{->}@/^4pc/@<-2ex>[dd]&\\
& & & &  T^2_{\gt{X}\mtn X} \ar@{:>}[dlll]\ar@{=>}[drrr]& & & & & \\
& T^2_{ X} & & & & & & 0 \quad &}$
 \caption{The cotangent braid for the normalization $n:\gt{X} \mtn
 X$}\label{dig2}
 \end{figure}
The map $\varphi^{\star}: T^0_Y \sr T^1_{\gt{X}\tru \gt{X} \mtn
Y/Y}\cong T^0_Y(\ohoa_{\gt{X}})$ and $\dot{\varphi}: T^0_{\gt{X}} \sr
T^0_Y(\ohoa_{\gt{X}})$ in Figure \ref{dig1} can be made explicit by
using the isomorphisms in Lemma \ref{bd1}.6.

$\varphi^{\star}: \bigoplus\limits_{j=1}^n
\c\{\gd{x}\}\tfrac{\partial}{\partial x_j} \sr
\bigoplus\limits_{j=1}^n \ohoa_{\gt{X}}\tfrac{\partial}{\partial
x_j}$ is componentwise  the structure map
$$g\tfrac{\partial}{\partial x_j} \longmapsto
\varphi^\sharp(g)\tfrac{\partial}{\partial x_j}.
$$
 $\dot{\varphi} :
\bigoplus\limits_{i=1}^r \ohoa_{\gt{X}}\tfrac{\partial}{\partial t_i}
\sr \bigoplus\limits_{j=1}^n \ohoa_{\gt{X}}\tfrac{\partial}{\partial
x_j}$ is given by the tangent map of $\varphi$,
$$
h(t_i).\tfrac{\partial}{\partial t_i} \longmapsto \sum_{j=1}^n
h(t_i).{\dot{\varphi}}^{(j)}_i (t_i)\tfrac{\partial}{\partial x_j},
$$  where $\dot{\varphi}^{(j)}_i(t_i) = \tfrac{\partial
\varphi^{(j)}_i}{\partial t_i}$ and
 ${\dot{\varphi}}^{(j)}  = (\dot{\varphi}^{(j)}_i)_{i=1}^r \in \ohoa_{\gt{X}}$. In particular, we have
\begin{equation}\label{pt3}
\varphi^{\star}(T^0_Y) = \bigoplus\limits_{j=1}^n \ohoa_X
\tfrac{\partial}{\partial x_j},
\end{equation}
\begin{equation}\label{pt4}
 \dot{\varphi}(T^0_{\gt{X}}) = \ohoa_{\gt{X}}\bigg(
\sum_{j=1}^n \dot{{\varphi}}^{(j)}  \tfrac{\partial}{\partial
x_j}\bigg).
\end{equation}

Since for any $\ohoa_X$-module $M$, $T^0_X(M) \subset T^0_Y(j_*M),$
we see that, in Figure \ref{dig2}, the maps
$$ \bigoplus\limits_{j=1}^n \ohoa_X\tfrac{\partial}{\partial
x_j} \supset  T^0_X \overset{n^\star}{\sr} T^1_{\gt{X}\tru \gt{X}
\mtn X/X}\cong T^0_X(\ohoa_{\gt{X}})\subset
\bigoplus\limits_{j=1}^n\ohoa_{\gt{X}}\tfrac{\partial}{\partial
x_j},$$
$$\dot{n} : \bigoplus\limits_{i=1}^r\ohoa_{\gt{X}}\tfrac{\partial}{ \partial t_i} \cong
T^0_{\gt{X}} \sr T^1_{\gt{X}\tru \gt{X} \mtn X/X}\subset
\bigoplus\limits_{j=1}^n\ohoa_{\gt{X}}\tfrac{\partial}{\partial
x_j}$$ are given by the same formulas as $\varphi^\star$ and
$\dot{\varphi}.$
 In particular, we have
 \begin{equation}\label{pt5}
 \dot{n}(T^0_{\gt{X}}) =
\ohoa_{\gt{X}}\big(\sum_{j=1}^n
\dot{{\varphi}}^{(j)}\tfrac{\partial}{\partial x_j}\big).
\end{equation}

\section{Computing $T^1$ and $T^2$}

We keep the assumptions and notations of the preceding section. We
shall explicitly compute $T^1$ and $T^2$ for the deformation functors
defined there. The case of plane curve singularities, which is
technically  much easier, was treated in \cite[Chapter II.2.4]{GLS}
and in positive characteristic in \cite{CGL}.

\thr \label{md1}  Let $I = \seq{f_1,\cdots, f_k} \subset \ohoa_Y =
\c\{x_1,\cdots,x_n\}$ be the ideal defining a reduced curve
singularity $X \subset Y = (\c^n,0)$. The following holds:
\begin{enumerate}
\item[1.]
\begin{itemize}
\item[(i)] $T^1_{\gt{X}\tru Y} \cong \Coker(T^0_Y \sr T^1_{\gt{X}\tru
\gt{X} \mtn Y/Y}) \cong \bigoplus\limits_{j=1}^n
(\ohoa_{\gt{X}}/\ohoa_X)\tfrac{\partial}{\partial x_j}$\\
 is a $\c$-vector space of dimension $n\delta$.

 \item[(ii)] $T^2_{\gt{X}\tru Y} = 0.$
 \end{itemize}
 \item[2.]
 \begin{itemize}
 \item[(i)] $T^1_{\gt{X}/ Y} \cong \Coker(\dot{\varphi}) \cong
\bigg(\bigoplus\limits_{j=1}^n
\ohoa_{\gt{X}}\tfrac{\partial}{\partial x_j} \bigg)\bigg/
\ohoa_{\gt{X}}\bigg(\sum_{j=1}^n
\dot{{\varphi}}^{(j)}\tfrac{\partial}{\partial x_j}  \bigg)$\\ is an
$ \ohoa_{\gt{X}}$-module of rank $n-1$.
 \item[(ii)] $T^2_{\gt{X}/Y} = 0.$
 \end{itemize}
 \item[3.]
\begin{itemize}
\item[(i)] $T^1_{\gt{X}\mtn Y}  \cong \Coker(T^0_{\gt{X}} \mtn
T^1_{\gt{X}\tru Y})$

 \hspace{0.98cm} $\cong \bigg(\bigoplus\limits_{j=1}^n
\big(\ohoa_{\gt{X}}/\ohoa_X\big)\tfrac{\partial}{\partial x_j}
\bigg)\bigg/ \big(\ohoa_{\gt{X}}/\ohoa_X\big)\bigg(\sum_{j=1}^n
\dot{{\varphi}}^{(j)}\tfrac{\partial}{\partial x_j}  \bigg)$\\
 is a finite
dimensional vector space of dimension
$$ \dim_\c T^1_{\gt{X}\mtn Y} = n\delta - m_1 = (n-3)\delta + e.$$
\item[(ii)] $T^2_{\gt{X}\mtn Y} = 0$.
\item[(iii)] If $X$ is smoothable and  unobstructed\footnotemark, $\dim_\c
T^1_{\gt{X}\mtn Y} = (n-3)\delta + \tau$. \footnotetext{$X$ is called
unobstructed if the base space of the semi-universal deformation of
$X$ is smooth. Complete intersections, 3-space curves, and Gorenstein
4-space curves $X$  are known to be unobstructed (even $T^2_X = 0)$
and smoothable.}
\item[(iv)]
 If $X$ is a complete
intersection then $\dim_\c T^1_{\gt{X}\mtn Y} = (n-3)\delta + \tau$;
if $X$ is a plane curve, $\dim_\c T^1_{\gt{X}\mtn Y} = \tau -
\delta$;  if $X$ is a 3-space curve,  $\dim_\c T^1_{\gt{X}\mtn Y} =
\tau;$ if $X$ is a $4$-space Gorenstein curve, $\dim_\c
T^1_{\gt{X}\mtn Y} = \tau + \delta.$
\end{itemize}
\item[4.]
\begin{itemize}
\item[(i)]   $T^1_X \cong \Coker \big(d^*:
\Hom_{\ohoa_X}(\Omega_Y^1\otimes\ohoa_X, \ohoa_X) \sr
\Hom_{\ohoa_X}(I/I^2, \ohoa_X)\big)$\\ is a finite dimensional vector
space of dimension $\tau$, where $$d^*\big(\sum_{j=1}^n a_j
\tfrac{\partial }{\partial x_j} \big)(f_i) = \sum_{j=1}^n a_j
\tfrac{\partial f_i}{\partial x_j}, \quad i=1,\cdots, k.$$

\item[(ii)] If $X$ is a complete intersection with $I = \seq{f_1,\cdots,f_{n-1}}$ then
 $$T^1_X \cong \ohoa_X^{n-1}\big/J(f_1,\cdots,f_{n-1})\ohoa_{X}^{n},$$
 where
$J(f_1,\cdots,f_{n-1})$ is the Jacobian matrix of $(f_1,\cdots,
f_{n-1})$.
\item[(iii)] If $X$ is a complete intersection, or if
$n \leq 3$, or if $X$ is Gorenstein and $n=4$, then $T^2_X = 0$.
\end{itemize}
\item[5.]
\begin{itemize}
\item [(i)]  $T^1_{\gt{X}\mtn X} \cong \Coker\big( T^0_{\gt{X}}\sr
T^1_{\gt{X}\tru X} \big)$  is a finite dimensional vector
 space,\\
$ \dim_\c T^1_{\gt{X}\mtn X}  = \dim_\c T^1_{\gt{X}\tru X} - m_1$\\
 $= \dim_\c T^1_{\gt{X}\mtn Y}  +  \dim_\c
 T^1_X(\ohoa_{\gt{X}}/\ohoa_X) - \dim_\c
 \Hom_{\ohoa_X}\big(I/I^2, \ohoa_{\gt{X}}/\ohoa_X\big).$
\item[(ii)] $T^2_{\gt{X}\mtn X} \cong
T^2_{\gt{X}\tru X} \cong T^1_X\big( \ohoa_{\gt{X}}/\ohoa_X\big)$ is a
finite dimensional vector space.
 \item[(iii)]$\dim_\c
T^1_{\gt{X}\mtn X} - \dim_\c T^2_{\gt{X}\mtn X} =  n\delta - m_1 -
\dim_\c \Hom_{\ohoa_X}\big(I/I^2, \ohoa_{\gt{X}}/\ohoa_X\big).$
\item[(iv)] If $X$ is a complete intersection then
 $$\dim_\c T^1_{\gt{X}\mtn X} - \dim_\c T^2_{\gt{X}\mtn X} =  \tau - 2\delta.$$
\item[(v)]  If $X$ is a plane curve singularity ($n=2)$ then\\
 \hspace{1cm}$T^1_{\gt{X}\mtn X} \cong T^1_{\gt{X}\mtn Y}  $ is of
dimension $\tau - \delta$; \\
\hspace{0.9cm} $T^2_{\gt{X}\mtn X}\cong \ohoa_{\gt{X}}/\ohoa_X $ is
of dimension $\delta$.
\end{itemize}
\item[6.]
\begin{itemize}
\item[ (i)]  $T^1_{\gt{X}\tru X} \cong T^0_X\big( \ohoa_{\gt{X}}/\ohoa_X\big) \cong \Hom_{\ohoa_X}\big(\Omega_X^1,  \ohoa_{\gt{X}}/\ohoa_X\big)$
  is a finite dimensional vector space of
dimension $\dim_\c T^1_{\gt{X}\mtn X} + m_1.$
\item[(ii)]  $T^2_{\gt{X}\tru X} \cong T^1_X\big( \ohoa_{\gt{X}}/\ohoa_X\big) $\\
$\cong  \Coker \Big(d^*: \Hom_{\ohoa_X}\big(\Omega_Y^1\otimes\ohoa_X,
\ohoa_{\gt{X}}/\ohoa_X \big)  \mtn \Hom_{\ohoa_X}\big(I/I^2,
\ohoa_{\gt{X}}/\ohoa_X \big)\Big) $
 is a finite dimensional vector space, and\\
$\dim_\c T^1_{\gt{X}\tru X} - \dim_\c T^2_{\gt{X}\tru X} = n\delta -
\dim_\c \Hom_{\ohoa_X}\big(I/I^2, \ohoa_{\gt{X}}/\ohoa_X \big).$
\item[(iii)] If $X$ is a complete intersection then $\dim_\c T^2_{\gt{X}\tru
X} \leq (n-1)\delta$ and
 $$\dim_\c T^1_{\gt{X}\tru X} - \dim_\c T^2_{\gt{X}\tru X} =
\delta.$$
\item[(iv)]  If $X$ is  a plane curve singularity then\\
\hspace{1cm} $T^1_{\gt{X}\tru X} \cong T^1_{\gt{X}\tru Y}$ is of
dimension $2\delta$; \\
\hspace{1cm} $T^2_{\gt{X}\tru X} \cong  \ohoa_{\gt{X}}\big/
\ohoa_X $ is of dimension $\delta$. \\
\end{itemize}
\item[7.]
\begin{itemize}
\item[(i)]  $ T^1_{\gt{X}/ X} \cong \Coker\big( T^0_{\gt{X}} \sr
T^0(\ohoa_{\gt{X}})\big)$

 \hspace{0.98cm}$\cong \ohoa_{\gt{X}} \ps{\sum_{j=1}^n
\dot{{\varphi}}^{(j)}\frac{\partial}{\partial
x_j}}{\gcd\big(\dot{\varphi}^{(1)},\cdots,\dot{\varphi}^{(n)}
\big)}\bigg/ \ohoa_{\gt{X}} \bigg( \sum_{j=1}^n
\dot{\varphi}^{(j)}\tfrac{\partial}{\partial x_j}\bigg)$

\hspace{0.9cm} $\cong \ps{\sum_{j=1}^n
\ohoa_{\gt{X}}.\big(t_i^{n^{(j)}_i -
mt_i}\big)_{i=1}^r\tfrac{\partial}{\partial x_j}}{\sum_{j=1}^n
\ohoa_{\gt{X}}.\big(t_i^{n^{(j)}_i -
1}\big)_{i=1}^r\tfrac{\partial}{\partial x_j}},$\\
 where $n^{(j)}_i =
\mbox{ord}_{t_i} \varphi^{(j)}_i(t_i), \quad j = 1,\cdots,n,
\quad i = 1,\cdots,r,$\\
 is a finite
dimensional vector space of dimension $mt - r$.
\item[(ii)] If $T^2_X = 0$ then
$$\dim_\c T^2_{\gt{X}/ X} = \tau + mt - r - \big( \dim_\c
T^1_{\gt{X}\mtn X} - \dim_\c T^2_{\gt{X}\mtn X} \big).$$
\item[ (iii)]  If $X$ is a complete intersection then
$$\dim_\c T^2_{\gt{X}/ X} = 2 \delta + mt - r = \mu + mt
- 1.$$
\end{itemize}
\item[8.]
\begin{itemize}
\item[ (i)]   $T^1_{\gt{X}\tru \gt{X} \mtn X/X} \cong
T^0(\ohoa_{\gt{X}}) \cong \ohoa_{\gt{X}} \ps{\sum_{j=1}^n
\dot{\varphi}^{(j)}\frac{\partial}{\partial
x_j}}{\gcd\big(\dot{\varphi}^{(1)},\cdots,\dot{\varphi}^{(n)}\big)}$\\
 is a free $\ohoa_{\gt{X}}$-module of rank $1$.
\item[ (ii)] $T^2_{\gt{X}/X}\cong T^2_{\gt{X}\tru \gt{X} \mtn X/X} \cong
T^1_X(\ohoa_{\gt{X}}) \cong $\\
 $ \cong \Coker \big(d^*:
\Hom_{\ohoa_X}(\Omega_Y^1\otimes\ohoa_X, \ohoa_{\gt{X}})
\sr \Hom_{\ohoa_X}(I/I^2, \ohoa_{\gt{X}})\big)$\\
  is a
finite dimensional vector space.
 \end{itemize}
\end{enumerate}
\ethr

\pf 1. (i)\quad  From the exact sequence $ \dashrightarrow $ in the
cotangent braid for the parametrization (Figure \ref{dig1}), we get
$$T^1_{\gt{X}\tru Y} \cong \Coker(\varphi^{\star}: T^0_Y \sr
T^1_{\gt{X}\tru \gt{X} \mtn Y/Y}) .$$  Then the formula follows from
the fact that $T^1_{\gt{X}\tru \gt{X} \mtn Y/Y} \cong
T^0_Y(\ohoa_{\gt{X}}) \cong \bigoplus\limits_{j=1}^n \ohoa_{\gt{X}}
\tfrac{\partial}{\partial x_j}$ and $(\ref{pt3})$. Since $\dim_\c
\ohoa_{ \gt{X}}/\ohoa_X = \delta$, we have
$\dim_\c T^1_{\gt{X}\tru Y} = n \delta.$\\
(ii) \quad $T^2_{\gt{X}\tru Y}$ appears in the same exact sequence of
complex vector spaces
$$
0 \dashrightarrow T^0_{\gt{X}\tru Y} \dashrightarrow T^0_{Y}
\dashrightarrow \cdots \dashrightarrow T^1_Y \dashrightarrow
T^2_{\gt{X}\tru \gt{X} \mtn Y/Y}\dashrightarrow T^2_{\gt{X}\tru
Y}\dashrightarrow T^2_Y \dashrightarrow\cdots $$ of the cotangent
braid for the parametrization. Noting that $T^2_Y = 0$ and
$$T^2_{\gt{X}\tru \gt{X} \mtn Y/Y} \cong T^1_{Y}(\ohoa_{\gt{X}})= 0,$$
since $Y
= (\c^n,0)$ is smooth, we get $T^2_{\gt{X}\tru Y} = 0.$\\
2.(i) and 2.(ii) follow in the same way from the exact sequence
$\longrightarrow $ in the cotangent braid for the parametrization
(Figure \ref{dig1})  and
$(\ref{pt4})$.\\
3.\quad  Applying the functor $\Hom_{\ohoa_{X}}\big( - ,
\ohoa_{\gt{X}}/\ohoa_X\big)$ to the defining exact sequence of
$\Omega^1_X$
\begin{equation}\label{pt2.6}
I/I^2 \overset{d}{\sr} \Omega_Y^1\otimes \ohoa_X \sr \Omega_X^1 \sr
0, \quad d(f_i) = \sum_{j=1}^n \tfrac{\partial f_i}{\partial x_j}
dx_j, \end{equation} noting that \big(by (\ref{pt2'}) and by
1.(i)\big)
$$ \Hom_{\ohoa_{X}}\big( \Omega_X^1 ,
\ohoa_{\gt{X}}/\ohoa_X\big)\cong T^0_X(\ohoa_{\gt{X}}/\ohoa_X) \cong
T^1_{\gt{X}\tru X},$$
$$\Hom_{\ohoa_{X}}\big(\Omega_Y^1\otimes \ohoa_X,
\ohoa_{\gt{X}}/\ohoa_X\big)\cong T^1_{\gt{X}\tru Y} \mbox{ and }
$$
$T^1_X(M) = \Coker\big( \Hom_{\ohoa_{X}}(\Omega_Y^1\otimes \ohoa_X ,
M \big) \sr \Hom_{\ohoa_{X}} ( I/I^2, M) \big)$  for each $\ohoa_X$-
module $M$, we get an exact sequence
$$T^1_{\gt{X}\tru X}\hookrightarrow T^1_{\gt{X}\tru Y}\overset{d^*}{
\longrightarrow} \Hom_{\ohoa_X}\big( I/I^2,
\ohoa_{\gt{X}}/\ohoa_X\big) \twoheadrightarrow
T^1_X(\ohoa_{\gt{X}}/\ohoa_X).$$ By combining this  with the exact
sequences $\Longrightarrow$ of the cotangent braids for $\gt{X} \sr
X$ and $\gt{X} \sr Y$ we obtain a  commutative diagram with exact
rows and columns as in Figure \ref{dig3}. The first vertical
isomorphism follows from diagram chasing.

\begin{figure} $$\xymatrix@C=15pt@R=9pt@M=5pt{
 0\ar[r] & T^0_{\gt{X}\mtn X}\ar[r]\ar[d]^{\cong}&
T^0_{\gt{X}}\ar[r]\ar @{=}[d]&
T^1_{\gt{X}\tru X} \ar[r]\ar @{^{(}->}[d] & T^1_{\gt{X}\mtn X}\ar[r] \ar @{^{(}->}[d] &0\\
0\ar[r] & T^0_{\gt{X}\mtn Y}/T^0_{\gt{X}\tru Y}\ar[r] &
T^0_{\gt{X}}\ar[r]& T^1_{\gt{X}\tru
Y}\ar[d]^{d^*}\ar[r] & T^1_{\gt{X}\mtn Y} \ar[r] &0\\
& & &{\phantom{A}}\save[]+<0cm,0cm>*{\Hom_{\ohoa_X}\big( I/I^2,
\ohoa_{\gt{X}}/\ohoa_X\big)} \ar@{->>}[d]
\restore \\
&&&{\phantom{A}}\save[]+<0cm,0cm>*{T^1_X(\ohoa_{\gt{X}}/\ohoa_X)\cong
T^2_{\gt{X}\tru X}} \restore }$$
 \caption{}\label{dig3}
\end{figure}
From the second row of the commutative diagram in Figure \ref{dig3}
we have
\begin{align}\label{pt3.2}
 T^1_{\gt{X}\mtn Y} &= \Coker(T^0_{\gt{X}} \sr
T^1_{\gt{X}\tru Y} )\nonumber \\
& = \Coker(T^0_{\gt{X}} \overset{\dot{\varphi}}{\sr} T^1_{\gt{X}\tru
\gt{X}\mtn Y/Y}\big/ \varphi^{\star}(T^0_{Y})).
\end{align}
Then the formula for
$T^1_{\gt{X} \mtn Y}$ follows from 2.(i),  (\ref{pt3}) and (\ref{pt4}). \\
Now we compute its dimension. By using the diagram in Figure
\ref{dig3}, and the fact that $T^0_{\gt{X}\mtn X} \cong T^0_{X}$
(see, Lemma \ref{bd1}.1) we get
$$\Im(T^0_{\gt{X}} \sr
T^1_{\gt{X}\tru Y})\cong T^0_{\gt{X}}\big/ T^0_{\gt{X}\mtn X} \cong
T^0_{\gt{X}}\big/ T^0_{X}. $$
 Hence, using 1.(i) we get 3.(i) from
\begin{align}
\dim_\c T^1_{\gt{X}\mtn Y}& =\dim_\c T^1_{\gt{X}\tru Y}-\dim_\c
\Im(T^0_{\gt{X}} \mtn T^1_{\gt{X}\tru Y}) \nonumber\\
 & = n\delta-\dim_\c
\big(T^0_{\gt{X}}\big/ T^0_{X}\big)= n\delta - m_1.\nonumber
\end{align}
 3.(ii) follows in the same way as 1.(ii) from the exact
sequence $\Longrightarrow$ in the cotangent braid for the
parametrization. \\
 3.(iii) follows from the fact that $e = \tau$  when $X$
is smoothable and unobstructed. \\
 Since complete intersections,  3-space curves and Gorenstein 4-space curves
  are smoothable and unobstructed, the proof for 3.(iv) follows. \\
4. is well-known (cf. \cite[Prop. II.1.25, II.1.29]{GLS}).\\
5.(i) \quad The formula  for $T^1_{\gt{X}\mtn X}$ and the first
formula for its dimension  follow from the exact sequence
$\Longrightarrow$ in the cotangent braid for the normalization
(Figure \ref{dig2}), noting that $\Im (T^0_{\gt{X}} \sr
T^1_{\gt{X}\tru X}) \cong T^0_{\gt{X}}\big/T^0_{\gt{X}\mtn X} \cong
T^0_{\gt{X}}\big/T^0_{X}$.

To show the second formula for its dimension, we note that the map
$$d^*:\Hom_{\ohoa_X}\big(\Omega_Y^1\otimes \ohoa_X,
\ohoa_{\gt{X}}/\ohoa_X\big) \sr \Hom_{\ohoa_X}\big(I/I^2,
\ohoa_{\gt{X}}/\ohoa_X\big)\big)$$ is given by the Jacobian matrix of
$(f_1,\cdots, f_k)$, $$d^*(\sum_{j=1}^n\limits a_j
\tfrac{\partial}{\partial x_j})(f_i) = \sum_{j=1}^n\limits
a_j\tfrac{\partial f_i}{\partial x_j}, a_j \in
\ohoa_{\gt{X}}/\ohoa_X,
 \quad i=1,\cdots, k.$$

The kernel of this map is $T^0_{X}\big(\ohoa_{\gt{X}}/\ohoa_X\big)$
by (\ref{pt2.6}). Hence, the dimension of $\Im(d^*)$ is $\dim_\c
T^1_{\gt{X}\tru Y} - \dim_\c T^0_{X}\big(\ohoa_{\gt{X}}/\ohoa_X\big)
= \dim_\c T^1_{\gt{X}\tru Y} -  \dim_\c T^1_{\gt{X}\tru X}$.

Moreover, the third column of the  diagram in Figure \ref{dig3} gives
us\\
 $T^2_{\gt{X}\tru X}  \cong
T^1_{X}\big(\ohoa_{\gt{X}}/\ohoa_X\big)$\\
 $\cong
\Coker(d^*:\Hom_{\ohoa_X}\big(\Omega_Y^1\otimes \ohoa_X,
\ohoa_{\gt{X}}/\ohoa_X\big) \sr \Hom_{\ohoa_X}\big(I/I^2,
\ohoa_{\gt{X}}/\ohoa_X\big)\big)$\\
 which also shows the first part of
6.(ii). Combining this  with 1.(i), 3.(i)  gives the second
formula for dimension of $T^1_{\gt{X}\mtn X}$. \\
5.(ii) follows from the exact sequence $\Longrightarrow$ in Figure
\ref{dig2}.\\
5.(iii) is a consequence of 5.(i) and 5.(ii). \\
In order to  prove 5.(iv), we use the fact that, when $X$ is a
complete intersection, $I/I^2$ is free of rank $n-1$. Then
$$\dim_\c\Hom_{\ohoa_X}\big(I/I^2, \ohoa_{\gt{X}}/\ohoa_X\big) =
(n-1)\delta.$$
 This implies, together with 3.(iii) and 5.(iii), the
statement of 5.(iv).\\
For proof of 5.(v), we notice from above that
$$T^1_{\gt{X}\mtn X} \cong \Coker(T^0_{\gt{X}}\sr T^1_{\gt{X}\tru
X}),$$
$$T^1_{\gt{X}\mtn Y} \cong \Coker(T^0_{\gt{X}}\sr T^1_{\gt{X}\tru
Y}).$$ Hence, it  suffices to show 6.(iv).\\
6.  The formula for $T^1_{\gt{X}\tru X}$ is obtained from
(\ref{pt2'}) and its dimension follows from 5.(i). 6.(ii) follows
from 6.(i),
the proof of 5.(i), and  5.(ii),  5.(iii).\\
If  $X$ is a complete intersection, $\dim_\c
\Hom_{\ohoa_X}\big(I/I^2, \ohoa_{\gt{X}}/\ohoa_X\big) = (n-1)\delta.$
Hence 6.(iii) follows from 6.(ii) and the surjection in the third
column
of the diagram in Figure \ref{dig3}.\\
Now we shall prove 6.(iv). We know that $d^*$ (in Figure \ref{dig3})
is induced by the Jacobian matrix of the defining function $f$ of
$X$. It is known (cf. \cite[Lemma II.2.31, p.316]{GLS})\footnotemark
~ \footnotetext{Another proof of this fact which does not use
deformation theory (but uses local duality)
 follows from applying $\wedge \mbox{d} f$ to the inclusion
 $n_*\Omega^1_{\gt{X}} \subset \omega^R_X$ which
 results in $\seq{\partial f /\partial x, \partial f/\partial y}
  \ohoa_{\gt{X}}\subset \ohoa_X$, where $\omega^R_X$
 are  Rosenlicht's regular  differential forms on $X$
  (cf. \cite[proof of Lemma II.2.32, p.317]{GLS}).}
that the Jacobian ideal of $f$ is contained in  the conductor
$\Ann_{\ohoa_X}\big( \ohoa_{\gt{X}}/\ohoa_X\big)$. By the formula for
$d^*$ in the proof of 5.(i),  $d^*$ is a zero map. This implies,
together with the third column of the diagram in Figure \ref{dig3},
the first isomorphism in 6.(iv) and
$$T^2_{\gt{X}\tru X} \cong T^1_{X}\big(
\ohoa_{\gt{X}}/\ohoa_X\big)\cong \Hom_{\ohoa_X}\big(I/I^2,
\ohoa_{\gt{X}}/\ohoa_X\big).$$ Since $I/I^2 \cong \ohoa_X f$, we get
$$T^2_{\gt{X}\tru X} \cong T^1_X \big(
\ohoa_{\gt{X}}/\ohoa_X\big) \cong \ohoa_{\gt{X}}/\ohoa_X$$ which has
$\c$-dimension $\delta$.\\
 7.(i) \quad  From
the exact sequence $\longrightarrow$ in the cotangent braid for the
normalization, we get
$$T^1_{\gt{X}/X} \cong \Coker\big(\dot{n}: T^0_{\gt{X}} \sr T^1_{\gt{X}\tru \gt{X}\mtn
X/X}\cong T^0_X(\ohoa_{\gt{X}})\big).$$ Consider the map
$$d^*: \Hom_{\ohoa_X}\big(\Omega_Y^1\otimes \ohoa_X,
\ohoa_{\gt{X}}\big) \sr  \Hom_{\ohoa_X}\big(I/I^2,
\ohoa_{\gt{X}}\big)$$ defined by
$$d^*\big(\sum_{j=1}^na_j \tfrac{\partial }{\partial
x_j}\big)(f_i) = \sum_{j=1}^na_j \tfrac{\partial f_i}{\partial x_j},
\quad a_j \in \ohoa_{\gt{X}},  \quad i = 1,\cdots, k.$$
 Note that
$T^0_X(\ohoa_{\gt{X}})$ is the kernel of this map \big(by applying
$\Hom_{\ohoa_X} (-, \ohoa_{\gt{X}})$ to  (\ref{pt2.6})\big) which is
torsion free, hence free, $\ohoa_{\gt{X}}-$module of rank 1. By the
chain rule, $\sum_{j=1}^n\tfrac{\partial f_i}{\partial
x_j}\dot{\varphi}^{(j)} = 0, \vm  \quad i =1,\cdots, k,$ hence $\big(
{\dot{\varphi}}^{(1)},\cdots, {\dot{\varphi}}^{(n)}\big) \in
\ker(d^*) = T^0_X(\ohoa_{\gt{X}})$, which is a non-zerodivisor (in
characteristic $0$). Therefore, the latter module is generated by
$$\big(
\dot{\varphi}^{(1)},\cdots, \dot{\varphi}^{(n)}\big)/ \gcd\big(
\dot{\varphi}^{(1)},\cdots, \dot{\varphi}^{(n)}\big).$$
  This implies,
together with $(\ref{pt5})$, the formula for $T^1_{\gt{X}/X}$. Since
$$\ps{\sum_{j=1}^n \ohoa_{\gt{X}}.\big(t_i^{n^{(j)}_i -
mt_i}\big)_{i=1}^r\tfrac{\partial}{\partial x_j}}{\sum_{j=1}^n
\ohoa_{\gt{X}}.\big(t_i^{n^{(j)}_i -
1}\big)_{i=1}^r\tfrac{\partial}{\partial x_j}} \cong
\bigoplus_{i=1}^r \c\{t_i\}/\seq{t_i^{mt_i - 1}}$$ which is of
$\c-$dimension $mt - r$, we get  the dimension of $T^1_{\gt{X}/X}$.
\\
(ii)\quad If $T^2_X = 0$, the exact sequence $\xymatrix{
\ar@{:>}[r]&}$ in Figure \ref{dig2} reads
$$ 0 \longrightarrow T^1_{\gt{X}/X} \longrightarrow T^1_{\gt{X}\mtn
X}\sr T^1_X \sr T^2_{\gt{X}/X} \sr T^2_{\gt{X}\mtn X} \sr 0.$$ By
taking the alternating sum of dimensions in this exact sequence,
together with the formula for the  dimension of $T^1_{\gt{X}/X}$ from
7.(i), we obtain 7.(ii).\\
Since $T^2_X = 0$ when $X$ is a complete intersection, 7.(iii)
follows from 7.(ii) and 5.(iv). \\
8.(i) is already shown in the proof of 7.(i).\\
The first isomorphism in 8.(ii) follows from the exact sequence
$\longrightarrow$ in Figure \ref{dig2}. The second isomorphism is
obvious by (\ref{pt1}). Finally, the third isomorphism follows from
the above statement about $T^1_X(M)$. \epf

\begin{remark} \label{rm3.2}\rm
  As a corollary to the dimension statements of Theorem 1, we get
interesting consequences concerning the semi-universal deformation of
the corresponding functors, either in the category of complex
analytic germs or in the category of complete analytic local rings
for algebraically closed fields $\k$ of characteristic $0$. We refer
to \cite{Bu} in the formal case and to \cite{Pal} and \cite{Fle} in
the complex analytic case. We formulate these consequences for both
cases in the geometric language where in the formal case the spaces
are to be understood as (formal) spectra of the complete local rings.

Recall that $\dim_\k T^1 < \vc$ for one of the deformation functors
of Section 2 holds if and only if it has a semi-universal deformation
(complex analytic respectively formal) and that $T^1$ is the Zariski
tangent space of the corresponding base space $B$ (also called
\emph{local moduli space}). $B$ is smooth in case $T^2 = 0$. More
generally we have the following estimate for the dimension of $B$,
$$\dim_\k T^1 - \dim_\k T^2 \leq \dim B \leq
\dim_\k T^1.$$

Furthermore, if a natural transformation of deformation functors
which have a versal deformation, induces a surjection on $T^1$ and an
injection on $T^2$, then the induced morphism on the semi-universal
base spaces is smooth (cf. \cite[Section II.1.3]{GLS} resp.
\cite{Bu}). Hence all eight deformation functors associated to
$\gt{X} \mtn X \mtn Y$ considered in Theorem \ref{md1} have a
semi-universal deformation except $\gd{\Def}_{\gt{X}/ Y}$ and
$\gd{\Def}_{\gt{X}\tru \gt{X} \mtn X/X}$, and we obtain exact
formulas resp. estimates for the dimension of the semi-universal base
spaces.

In particular, from Theorem \ref{md1}, 3.(i) and (ii) we get that the
parametrization $\varphi: \gt{X} \mtn Y$ of $X$ has a semi-universal
deformation with smooth base space $B$ of dimension
$$ \dim B = n\delta - m_1 = (n-3)\delta + e.$$

Note that $B$ is also the base space of the semi-universal unfolding
of $\varphi$ with respect to (extended) left-right equivalence (cf.
Section 5).

 If we consider $\varphi$ only with respect to left
equivalence, the base space of its semi-universal unfolding coincides
with the semi-universal base space of the functor
$\gd{\Def}_{\gt{X}\tru Y}$ and hence it is smooth of dimension
$n\delta$ according to Theorem \ref{md1}, 1.(i) and (ii).
\end{remark}
\section{Applications to $\hoa{A}$-equivalence of parametrized curve singularities}

Let $\hoa{A}$  denote the Mather group of left-right equivalence of
morphisms between manifolds (cf. \cite{Ma}). We consider
(multi-)germs of analytic maps
$$\varphi: (\c, S) \sr (\c^n,0), \quad S = \{z_1,\cdots,z_r\},$$
 where $(\c,S)$ is the multigerm $\coprod_{i=1}^r (\c,z_i)$. After
 choosing local coordinates $s_i$ of $(\c,z_i)$ with center $0$ and
 $x_1,\cdots, x_n$ of $(\c^n,0)$, $\varphi$ is given as
$$\varphi = \big(\varphi_i\big)_{i=1}^r, \quad  \varphi_i = \big(\varphi_i^{(1)},\cdots,
\varphi_i^{(n)}\big): (\c,0) \sr (\c^n,0), $$ where the $\varphi_i$
are holomorphic map germs. Then $\varphi^\sharp : \ohoa_{\c^n,0} \sr
\ohoa_{\c,S}$ is given, as in Section 3, by $\varphi^\sharp(x_j) =
\varphi^{(j)} = (\varphi^{(j)}_1,\cdots, \varphi^{(j)}_r).$ Applying
$\tfrac{d}{d s_i}$ to the $i$-th component, we get an element
$$\dot{\varphi}^{(j)} = \big(\dot{\varphi}^{(j)}_1,\cdots,\dot{\varphi}^{(j)}_r
\big) \in \bigoplus\limits_{i=1}^r \c\{s_i\} = \ohoa_{\c,S}.$$
 We denote by $$\varphi: \nga{X} = \coprod_{i=1}^r \nga{X_i} \sr Y$$ a
representative of the multigerm of $\varphi$, where $\nga{X_i}$ resp.
$Y$ are (sufficiently small) open neighborhoods of $0 \in \c$ resp.
$0 \in \c^n.$ Moreover, we write also $(\nga{X},\nga{0})$ instead of
$(\c,S)$. \\
We do not assume, as in Section 3, that $\varphi$ factors through the
normalization of a reduced curve singularity. However, as we show,
this condition is necessary and sufficient for  $\varphi$ being
finitely $\hoa{A}$-determined (which clarifies the various sufficient
conditions in the literature). The aim is to apply the results of the
previous sections, obtained by deformation theory, to determine
formulas resp. estimates for the codimension of $\varphi$ (in the
sense of Mather and Wall) when $\varphi$ is finitely
$\hoa{A}$-determined.

\quad
Let us first establish a dictionary from  Mather's to our notations:\\
$\theta(\varphi) = T^0_{Y,0}(\ohoa_{\nga{X},\nga{0}}) \cong
\bigoplus\limits_{j=1}^n
\ohoa_{\nga{X},\nga{0}}\tfrac{\partial}{\partial x_j}$, germs
 of vector fields $\sigma: (\c,S) \sr T\c^n$ along  the parametrization
 $\varphi$.\\
$ \theta_{\c, S} =  T^0_{\nga{X},\nga{0}} \cong
\bigoplus\limits_{i=1}^r\c\{s_i\}\tfrac{\partial}{\partial s_i}, $
germs of vector fields along the identity in $(\c, S)$.\\
 $\theta_{\c^n, 0} = T^0_{Y,0} \cong \bigoplus\limits_{j=1}^n \ohoa_{Y,0}\tfrac{\partial}{\partial x_j},$
 germs of vector fields along the identity in $(\c^n,0).$\\
 Moreover, we have the linear maps $t\varphi$ resp. $\omega\varphi$
 of Mather, coinciding with our $\dot{\varphi}$ resp.
 $\varphi^{\star}$  in Section 3:
$$\dot{\varphi} = t\varphi: \theta_{\c, S} \sr \theta(\varphi),\quad
t\varphi(\xi) = d\varphi \circ \xi,$$
$$\varphi^\star = \omega \varphi: \theta_{\c^n, 0} \sr \theta(\varphi), \quad
\omega \varphi(\eta) = \eta\circ \varphi.$$ The \emph{extended
tangent space } to the orbit of the multigerm $\varphi: (\c,S) \sr
(\c^n,0)$ by the action of the  group $\hoa{A}$ is defined as (cf.
\cite[Part I.1]{Wa})
 $$T_e \hoa{A}\varphi = t\varphi(\theta_{\c, S}) +
\omega\varphi(\theta_{\c^n, 0}),$$ and
$$ d_e(\varphi, \hoa{A}):= \dim_\c \big(\theta(\varphi)/T_e\hoa{A}\varphi\big)$$ the
\emph{$\hoa{A}_e$-codimension }of $\varphi$. By \cite[Theorem
1.2]{Wa} we have that  $\varphi$ is finitely $\hoa{A}$-determined if
and only if $d_e(\varphi,\hoa{A}) < \vc$.

The relation between deformations of $\varphi: (\c,S) \mtn (\c^n,0)$
as in the previous sections and $\hoa{A}$-equivalence of $\varphi$ is
as follows. \\
Elements of the extended tangent space $T_e\hoa{A}\varphi\subset
\theta(\varphi)$ correspond to infinitesimal deformations of
$\varphi$ as elements of $\Def_{(\c,S)\mtn (\c^n,0)}(T_\epsilon)$
which are locally trivial. Dividing out the locally trivial
deformations from the latter we get the vector space
$\gd{\Def}_{(\c,S)\mtn (\c^n,0)}(T_\epsilon) = T^1_{(\c,S)\mtn
(\c^n,0)}.$ By passing to the quotients we get in fact an isomorphism
of vector spaces,
$$\theta(\varphi)/T_e\hoa{A}\varphi \cong T^1_{(\c,S)\mtn
(\c^n,0)}$$ and, in particular,
\begin{equation}\label{pt41}
d_e(\varphi, \hoa{A}) = \dim_\c T^1_{(\c,S)\mtn (\c^n,0)}.
\end{equation}
This follows from (\ref{pt3.2}), (\ref{pt1}), Lemma \ref{bd1}.6 and
the dictionary above.

 Using  Lemma \ref{bd1}.6, and  formulas (\ref{pt3}),
(\ref{pt4}) which hold for arbitrary $\varphi$, we obtain
 \lm  Let $\varphi: (\c,S) \sr (\c^n,0)$ be a
 complex analytic multigerm. Then
$\quad \theta(\varphi) \cong \bigoplus\limits_{j=1}^n \ohoa_{\c,S}
\tfrac{\partial}{\partial x_j}, \quad  t\varphi(\theta_{\c,S}) \cong
\ohoa_{\c,S}.\big(\sum_{j=1}^n
\dot{\varphi}^{(j)}\tfrac{\partial}{\partial x_j}\big),$\\
$\omega\varphi(\theta_{\c^n,0}) \cong \bigoplus\limits_{j=1}^n
(\varphi^\sharp\ohoa_{\c^n,0})\tfrac{\partial}{\partial x_j},$ and
\begin{equation}\label{pt4.1}
 d_e(\varphi,\hoa{A}) = \dim_\c
\ps{\bigoplus\limits_{j=1}^n \ohoa_{\c,S} \tfrac{\partial}{\partial
x_j}}{\ohoa_{\c,S}.\big(\sum_{j=1}^n
\dot{\varphi}^{(j)}\tfrac{\partial}{\partial x_j}\big) +
\bigoplus\limits_{j=1}^n
(\varphi^\sharp\ohoa_{\c^n,0})\tfrac{\partial}{\partial x_j}}.
\end{equation}
 \elm

\begin{remark}\rm\begin{itemize}
 \item[(1)] To determine when $d_e(\varphi, \hoa{A})  < \vc,$  note that
$\varphi_i: (\c,0) \sr (\c^n,0)$  must not be constant for any $i$.
Because otherwise $\varphi_i = \dot{\varphi}_i$ is identically zero
and hence the $i-$th component of the denominator in (\ref{pt4.1}) is
$0$, implying $d_e(\varphi,\hoa{A}) = \vc.$

 \item[(2)] $\varphi_i$ being non-constant is equivalent to $\dim
\varphi_i^{-1}(0) = 0$ and hence to $\varphi_i$ being finite (cf.
\cite[Theorem 1.70]{GLS}). Furthermore, this is equivalent to
$(X_i,0):= \varphi_i(\c,0) \subset (\c^n,0)$ being the germ of curve
singularity (cf. \cite[Corollary 1.68]{GLS}). The map $\varphi$ is
finite iff $\varphi_i$ is finite for all $i$.

\item[(3)]  Consider the case that $\varphi_i$ is an isomorhism onto
the same smooth image for all i. By choosing appropriate coordinates,
we may assume that $\varphi_i(s_i)=(s_i,0,...,0)$ for all i. Then,
using (\ref{pt4.1}), we see that  $d_e(\varphi,\hoa{A}) = 0$ if $r=1$
but that $d_e(\varphi,\hoa{A}) = \vc$ if $r>1$.

\end{itemize}
\end{remark}

If $\varphi_i$ is finite we equip $(X_i,0)$ with its reduced
structure, that is, the ideal $I_i\subset \ohoa_{\c^n,0}$ of
$(X_i,0)$ is the vanishing ideal and call $\varphi_i$ a
\emph{parametrization} of the branch $(X_i,0)$. It is called
\emph{primitive} if for any other parametrization $\nga{\varphi_i}$
of $(X_i,0)$, satisfying $\varphi_i(s_i) =
\nga{\varphi_i}(\nga{s_i}^k)$ for some regular parameter $\nga{s_i} =
s_i.u_i, u_i \in \c\{s_i\}^*,$ we have $k=1$.

If the map $\varphi = (\varphi_i)_{i=1}^r$ is finite we set $(X,0) =
\bigcup_{i=1}^r (X_i,0)$ and call $\varphi$ a \emph{parametrization
of} $(X,0)$. A parametrization is called \emph{primitive} if
$\varphi_i$ is primitive for each  $i = 1,\cdots, r$ and if
$\varphi_i(\c,0) \not = \varphi_j(\c,0)$ for all $i \not = j$.

We denote by $X_i$ resp. $X = \varphi(\nga{X})$ sufficiently small
representatives of $(X_i,0)$ resp. $(X,0)$ and by $n: \nga{X} \mtn X$
the induced map.

To obtain necessary and sufficient conditions for $\varphi$ to be
finitely $\hoa{A}$-determined, we  consider the following sheaves on
$X$:
$$\Theta_\varphi := \bigoplus\limits_{j=1}^n \big(\varphi_*\ohoa_{\nga{X}}\big)\tfrac{\partial}{\partial x_j},
\quad \hoa{T}_\varphi := \varphi_*\ohoa_{\nga{X}}\big(\sum_{j=1}^n
\dot{\varphi}^{(j)} \tfrac{\partial}{\partial x_j} \big),$$
$$\hoa{W}_\varphi := \bigoplus\limits_{j=1}^n \big(\varphi^\sharp \ohoa_{Y}\big)\tfrac{\partial}{\partial x_j},
\quad \hoa{F} := \ps{\Theta_\varphi }{\hoa{T}_\varphi +
\hoa{W}_\varphi}.$$

If $\varphi$ is finite, $\hoa{F}$ is coherent (cf. \cite[Theorem
I.1.67]{GLS}) and we have $d_e(\varphi, \hoa{A}) = \dim_\c
\hoa{F}_0.$ Therefore, $d_e(\varphi,\hoa{A}) < \vc$ if and only if
$\hoa{F}$ is concentrated on $\{0\}$ (cf. \cite[Corollary
I.1.74]{GLS}).

\begin{proposition} \label{bd4.2} With the above notations let
$\varphi: (\nga{X},\nga{0}) = (\c,S) \sr (\c^n,0) = (Y,0)$, $ n \geq
2,$ be a (multi-) germ of a complex analytic morphism. The following
are equivalent: \\
(1) $d_e(\varphi,\hoa{A}) < \vc$ (i.e., $\varphi$ is finitely
$\hoa{A}$-determined). \\
(2) $\varphi$ is a primitive parametrization of $(X,0)$. \\
(3) $\varphi$ is finite and the sheaf $\hoa{F}$ is concentrated on $\{0\}$.\\
(4) $\varphi$ is finite and $n: (\nga{X},\nga{0}) \sr (X,0)$ is the
normalization of $(X,0)$. \\
(5) $\varphi$ is finite and $n: \nga{X}\tru \{\nga{0}\} \sr X \tru
\{0\}$ is bijective (i.e., $n$ is birational).
\end{proposition}

\pf If $\varphi$ is a parametrization of $(X,0)$ then the equivalence
of (4), (5) and (2) is a property of the normalization
(cf. \cite[Chapter I.1.9]{GLS}).\\
By Remark 2.(1) condition (1) implies that $\varphi$ is finite.
Therefore the equivalence of (1) and  (3) follows from the coherence
of $\hoa{F}$ as shown above.

It remains to show the equivalence $(2) \td (3)$. By Remark 2.(3),
applied to a non-singular point of $(X,0)$, $\hoa{F}$ is not
concentrated on $\{0\}$ if $\varphi_i(\c,0) = \varphi_j(\c,0)$ for
some $i \not = j$. Hence we have to show that $\hoa{F}_p \not= \{0\}$
for $p \in X_i\tru\{0\}$ if and only if $\varphi_i$ is not primitive
for some $i$. We may assume that $X$ is irreducible, and that
$\varphi_i = \varphi$. Then $\varphi$ factors as
$$\varphi: (\nga{X},0) = (\c,0) \overset{\nu}{\sr} (\c,0) =
(\gt{X},0) \overset{\gt{\varphi}}{\sr } (\c^n,0) = (Y,0)$$ with
$\gt{\varphi}$ primitive and $\nu(s) = s^k$ for some $k\geq 1$. Since
$\gt{\varphi}$ is an injective immersion outside $0$, we have
$\gt{\varphi}^{-1}(p) = \{q\}$ and
$$\Theta_{\varphi,p} \cong \bigoplus\limits_{j=1}^n \big(\nu_*\ohoa_{\nga{X}}\big)_q\tfrac{\partial}{\partial x_j}, $$
$$\hoa{T}_{\varphi,p} \cong (\nu_*\ohoa_{\nga{X}})_q\cdot\big(\sum_{j=1}^n \dot{\varphi}^{(j)} \tfrac{\partial}{\partial x_j}
\big), \quad \hoa{W}_{\varphi,p} \cong \bigoplus\limits_{j=1}^n
\big(\nu^\sharp \ohoa_{\gt{X},q}\big)\tfrac{\partial}{\partial
x_j}.$$ If $\nu^{-1}(q) = \{q_1,\cdots,q_k\}$ then $\nu:
(\nga{X},q_i)\sr (\gt{X},q)$ is an isomorphism for $i = 1,\cdots, k$
and hence $\nu^\sharp: \ohoa_{\gt{X},q} \sr \ohoa_{\nga{X},q_i}$ is
also an isomorphism, that is, $\nu^\sharp$ embeds  $\ohoa_{\gt{X},q}$
"diagonally" into $(\nu_*\ohoa_{\nga{X}})_q \cong
\bigoplus\limits_{i=1}^k \ohoa_{\nga{X},q_i}$. \\
It follows that $\Theta_{\varphi,p}$ is a free
$\ohoa_{\gt{X},q}$-module of rank $nk$ and $\hoa{W}_{\varphi,p}$ a
free $\ohoa_{\gt{X},q}$-module of rank $n$. \\
Since $\dot{\varphi}^{(j)}(s) =
\dot{\gt{\varphi}}^{(j)}(\nu(s)).\dot{\nu}(s)$ with $\dot{\nu}(q_i)
\not = 0$, $\dot{\nu}$ is a unit in $\ohoa_{\nga{X},q_i}$ and hence
$$\hoa{T}_{\varphi,p} \cong (\nu_*\ohoa_{\nga{X}})_q. \big(\sum_{j=1}^n \dot{\gt{\varphi}}^{(j)}\circ \nu
\tfrac{\partial}{\partial x_j} \big).$$ Since $\gt{\varphi}$ is an
isomorphism at $q$, $\sum_{j=1}^n \dot{\gt{\varphi}}^{(j)}\circ \nu
\tfrac{\partial}{\partial x_j}$ is a non-zero vector in
$\Theta_{\varphi,p}$ and therefore $\hoa{T}_{\varphi,p}$ is
$\ohoa_{\gt{X},q}$-free of rank $k$. Since
$$\hoa{T}_{\varphi,p} \cap \hoa{W}_{\varphi,p} \cong
\nu^\sharp\ohoa_{\gt{X},q}\seq{\dot{\gt{\varphi}}^{(j)}\circ \nu
\tfrac{\partial}{\partial x_j}}$$ is $\ohoa_{\gt{X},q}$-free of rank
1 we get that $\hoa{F}_p$ (which is free for $p\not = 0$ and
sufficiently close to $0$ by \cite[Theorem I.1.80]{GLS}) is of rank
$nk - n - k + 1 = (n-1)(k-1)$. Then $\hoa{F}_p \not= \{0\}$  if and
only if $k\geq 2$, that is, if and only if $\varphi$ is not
primitive. \epf

\thr  \label{dl4.1} Let $\varphi = (\varphi_i)_{i=1}^r : (\c,S) \sr
(\c^n,0)$ be a primitive parametrization of a reduced curve
singularity $(X,0)\subset (\c^n,0)$. Then
$$d_e(\varphi,\hoa{A}) =  n\delta
- m_1 = (n-3)\delta + e.$$
 Moreover we have the inequalities
$$(n-2)\delta \leq (n-2)\delta + t - 1 +  mt - r \leq n\delta - c + mt - r
\leq d_e(\varphi,\hoa{A}),$$
$$d_e(\varphi,\hoa{A}) \leq (n-1)\delta + \mu - c \leq n\delta - r <
n\delta.$$ \ethr Here and below $\delta,  \mu,  mt,  \tau,  $ etc.
are the invariants defined in Section 3; $c = \dim_\c
\big(\ohoa_{\c,S}/ \Ann(\ohoa_{\c,S}/\ohoa_{X,0}) \big)$ is the
degree of the conductor and $t$  is the Cohen-Macaulay type  of
$(X,0)$ (cf. \cite{Gr}).

\pf The formula for $d_e(\varphi,\hoa{A})$ is obtained from
(\ref{pt41}) and Theorem 1, 3.(i).  The inequalities follow from
\cite[Theorem 2.5]{Gr}.  \epf

\coro \label{hq4.2}
 \begin{itemize}
\item[(1)] If $(X,0)$ is smoothable and unobstructed then
$d_e(\varphi,\hoa{A}) = (n-3)\delta + \tau.$
\item[(2)] If $(X,0)$ is Gorenstein then $d_e(\varphi,\hoa{A}) \leq
(n-1)\delta - r + 1$ with equality if and only if $(X,0)$ is
quasihomogeneous.

\item[(3)] If $(X,0)$ is quasihomogeneous then $d_e(\varphi,\hoa{A}) =
(n-1)\delta - r + t.$

\item[(4)] $d_e(\varphi,\hoa{A}) = 0$ iff $(X,0)$ is smooth or an
ordinary node.
\end{itemize}

 \ecoro

\pf (1) follows since $e = \tau$ if $(X,0)$ is smoothable and
unobstructed (cf. \cite{De}, \cite{Gr}).\\
(2) By  \cite[Theorem 2.5 (2)]{Gr}, $e \leq \mu$ for $(X,0)$
Gorenstein. Hence $d_e(\varphi,\hoa{A}) = (n-3)\delta + e \leq
(n-3)\delta + \mu = (n-1)\delta - r + 1$ with equality if and only if
$(X,0)$ is quasihomogeneous by \cite[Satz 2.1]{GMP}. Statement
(3) follows from \cite[Theorem 2.5 (3)]{Gr}. \\
To see (4) let $(X,0)$ be singular. Theorem 2 implies
$d_e(\varphi,\hoa{A}) >0$ if $n>3$. If $n \leq 3$ then $(X,0)$ is
smoothable and unobstructed and we can apply (1). If $n = 3$ then
$d_e(\varphi,\hoa{A}) = \tau >0$ since $(X,0)$ is singular. If $n =
2$ then $d_e(\varphi,\hoa{A}) = \tau -\delta >0$ if and only if
$(X,0)$ is singular and not an ordinary node by \cite[Proposition and
Corollary]{Gr2}.  \epf

\begin{remark}\rm
\begin{itemize}
 \item[(1)] Whether the inequality $d_e(\varphi,\hoa{A}) \leq
 (n-1)\delta - r + t$ holds in general and whether  equality
 implies that $(X,0)$ is quasihomogeneous is an open problem. By  \cite[Remark 2.6 (2)]{Gr}
 these statements hold if $(X,0)$ is irreducible and if $t \leq 2$
 (note that $t = 1$ if and only if $(X,0)$ is Gorenstein).
 \item[(2)] The formulas of Corollary \ref{hq4.2}, (1) and (2) hold for complete intersections
 which
 are smoothable, unobstructed and Gorenstein. Corollary \ref{hq4.2} (2)
  was proved before for plane curve
 singularities  $(n=2)$ by Mond \cite[Theorem 2.3]{Mo}.
 \item[(3)] Corollary \ref{hq4.2} (3) generalizes \cite[Theorem 2]{HRR} where  the
 case of irreducible monomial curve singularities was considered.
 \item[(4)] Since $3$-space curve singularities are smoothable and
 unobstructed we get  $ d_e(\varphi,\hoa{A}) = \tau$ for $n=3$.
 \item[(5)] If $(X,0)$ is Gorenstein and $n=4$, then $(X,0)$ is  again
 smoothable and unobstructed, hence we get $d_e(\varphi,\hoa{A}) =
 \tau + \delta.$
 \item[(6)] In \cite[Proposition 4.52]{Wa} Wall showed (attributing the result to
 W.Bruce) that $d_e(\varphi,\hoa{L}) = \mu = 2\delta$ if $n=2$ and
 $(X,0)$ is irreducible, mentioning that it seems to be a curiosity
 that $d_e(\varphi,\hoa{L}) = d_e(f,\hoa{R})$ where $f$ is a function
 defining $(X,0)\subset (\c^2,0)$. Here $\hoa{L}$ resp. $\hoa{R}$
 denote the Mather groups of left resp. right equivalence. \\
 Indeed, $d_e(\varphi,\hoa{L})$ is in general not related to $\mu$ but to
 $\delta$ as the following result shows. Since left equivalence is the equivalence
 relation for $\gd{\Def}_{X\tru Y}$, we obtain from Theorem
 \ref{md1},
 1.(i) that   an arbitrary primitive parametrization is finitely
 $\hoa{L}$-determined and that
 $$d_e(\varphi,\hoa{L}) = \dim_\c T^1_{(\c,S)\tru (\c^n,0)} =
 n\delta.$$
\end{itemize}
\end{remark}

\vspace{0.5 cm}

{\bf Acknowledgement:} We like to thank the referee for his careful
reading of the manuscript and for his suggestions to section 5.



\begin{thebibliography}{}
\bibitem  [Bu]   {Bu} Buchweitz, R.O.:
\emph{Contributions \`{a} la th\'{e}orie de singularit\'{e}s.}
Th\`{e}se, Universit\'{e} Paris VII, 1981.
\bibitem  [CGL]   {CGL} Campillo, A.; Greuel, G.-M.; Lossen, C.:\emph{ Equisingular deformations of plane curves
 in arbitrary
characteristic.} Composito Math. \textbf{143}, 829-882 (2007).
\bibitem  [De]   {De} Deligne, P.: \emph{Intersection sur les surfaces reguliers.}
In. SGA VII 2, SLN \textbf{340} 1-38 (1973).
 \bibitem  [Fle] {Fle} Flenner, H.: \emph{\"{U}ber Deformationen holomorpher Abbildungen.} Habilitationsschrift,
Universit\"{a}t Osnabr\"{u}ck (1978).
\bibitem  [GLS]   {GLS} Greuel, G.-M.; Lossen, C.; Shustin, E.: \emph{Introduction to Singularites and
Deformations.} Springer (2007).
\bibitem  [GMP]   {GMP} Greuel, G.-M.; Martin, B.; Pfister, G.:
\emph{Numerische Charakterisierung Quasihomogener
Gorenstein-Kurvensingularit\"{a}ten.} Math. Nachr. \textbf{124},
123-131 (1985).
\bibitem  [Gr1]   {Gr} Greuel, G.-M.: \emph{On deformation of curves and a formula of Deligne.}
Algebraic geometry, Proc. int. Conf., La Rabida/Spain 1981,
Lect.Notes Math. \textbf{961}, 141-168 (1982).
\bibitem  [Gr2]   {Gr2} Greuel, G.-M.: \emph{A Remark on the Paper of
A. Tannenbaum.} Compositio Mat. \textbf{51}, 185-187 (1984).
\bibitem  [HRR]   {HRR} Hernandes, M.E.; Rodrigues Hernandes, M.E.; Ruas, M.A.S.:
\emph{$\hoa{A}_e$-codimension of germs of analytic  curves.} Manuscr.
Math. \textbf{124}, 237-246 (2007).
\bibitem  [Ma]   {Ma} Mather, J. N.: \emph{ Stability of $C^\infty$ mappings, III. Finitely determined
map-germs.} Publ. Math. IHES, \textbf{35 }, p. 127-156 (1968).
\bibitem  [Mo]   {Mo} Mond, D.: \emph{Looking at bent wires - $\hoa{A}_e$-codimension and the
vanishing topology of parametrized curve singularities.} Math. Proc.
Camb. Phil. Soc., \textbf{117}, 213-222 (1995).
 \bibitem [Pal]{Pal} Palamodov, V.P.: \emph{Deformations of complex spaces.} In: Several
 complex variables IV, Springer, Encyclopedia of math. sciences, Vol. \textbf{10} (1990).
 \bibitem  [Wa]   {Wa} Wall, C.T.C.: \emph{Finite determinacy of smooth
 map-germs.} Bull. London Math. Soc. (6) \textbf{13}, 481-539 (1981).
  \end{thebibliography}
\end{document}